\documentclass[12pt, a4paper]{article}
\usepackage{mymatrix}
\usepackage[cp1251]{inputenc}
\usepackage{amsmath}
\usepackage{amsthm}
\usepackage{euscript}
\usepackage{amssymb} \usepackage[dvips]{graphicx}
\begin{document}

\newcounter{bnomer} \newcounter{snomer}
\renewcommand{\thesnomer}{\thebnomer.\arabic{snomer}}
\renewcommand{\refname}{\begin{center}\large{\textbf{References}}\end{center}}

\newcommand{\sect}[1]{%
\setcounter{snomer}{0} \refstepcounter{bnomer}
\begin{center}\large{\textbf{\arabic{bnomer}. {#1}}}\end{center}}
\newcommand{\defi}[1]{%
\refstepcounter{snomer}
\par\textbf{Definition \arabic{bnomer}.\arabic{snomer}. }{#1}\par}
\newcommand{\theo}[1]{%
\refstepcounter{snomer}
\par\textbf{Theorem \arabic{bnomer}.\arabic{snomer}. }{#1} $\square$\par}
\newcommand{\theop}[2]{%
\refstepcounter{snomer}
\par\textbf{Theorem \arabic{bnomer}.\arabic{snomer}. }{#1}
\par\textbf{Proof}. {#2} $\square$\par}
\newcommand{\theosp}[2]{%
\refstepcounter{snomer}
\par\textbf{Theorem \arabic{bnomer}.\arabic{snomer}. }{#1}
\par\textbf{Sketch of the proof}. {#2} $\square$\par}
\newcommand{\exam}[1]{%
\refstepcounter{snomer}
\par\textbf{Example \arabic{bnomer}.\arabic{snomer}. }{#1}\par}
\newcommand{\deno}[1]{%
\refstepcounter{snomer}
\par\textbf{Notation \arabic{bnomer}.\arabic{snomer}. }{#1}\par}
\newcommand{\post}[1]{%
\refstepcounter{snomer}
\par\textbf{Proposition \arabic{bnomer}.\arabic{snomer}. }{#1} $\square$\par}
\newcommand{\postp}[2]{%
\refstepcounter{snomer}
\par\textbf{Proposition \arabic{bnomer}.\arabic{snomer}. }{#1}
\par\textbf{Proof}. {#2} $\square$\par}
\newcommand{\lemm}[1]{%
\refstepcounter{snomer}
\par\textbf{Lemma \arabic{bnomer}.\arabic{snomer}. }{#1} $\square$\par}
\newcommand{\lemmp}[2]{%
\refstepcounter{snomer}
\par\textbf{Lemma \arabic{bnomer}.\arabic{snomer}. }{#1}
\par\textbf{Proof}. {#2} $\square$\par}
\newcommand{\coro}[1]{%
\refstepcounter{snomer}
\par\textbf{Corollary \arabic{bnomer}.\arabic{snomer}. }{#1} $\square$\par}
\newcommand{\corop}[2]{%
\refstepcounter{snomer}
\par\textbf{Corollary \arabic{bnomer}.\arabic{snomer}. }{#1}
\par\textbf{Proof}. {#2} $\square$\par}
\newcommand{\nota}[1]{%
\refstepcounter{snomer}
\par\textbf{Remark \arabic{bnomer}.\arabic{snomer}. }{#1}\par}

\newcommand{\Ind}[3]{%
\mathrm{Ind}_{#1}^{#2}{#3}}
\newcommand{\Res}[3]{%
\mathrm{Res}_{#1}^{#2}{#3}}
\newcommand{\epsi}{\varepsilon}
\newcommand{\w}{\widetilde}
\newcommand{\Supp}[1]{%
\mathrm{Supp}\mathop{#1}}

\newcommand{\reg}{\mathrm{reg}}
\newcommand{\sreg}{\mathrm{sreg}}
\newcommand{\codim}{\mathrm{codim}\mathop}
\newcommand{\chara}{\mathrm{char}\,}
\newcommand{\vfi}{\varphi}
\newcommand\pho{\phantom{\otimes}}

\newcommand{\lee}{\leqslant}
\newcommand{\gee}{\geqslant}
\newcommand{\F}{\mathbb{F}}
\newcommand{\g}{\mathfrak{g}}
\newcommand{\p}{\mathfrak{p}}
\newcommand{\K}{\EuScript{K}}
\newcommand{\Ou}{\mathcal{O}}

\author{Mikhail V. Ignatyev\thanks{This research was supported by
Samara regional grant for students and young scientists no.~23E2.1D}}
\date{}
\title{Subregular characters of the unitriangular group \\over a finite field}
\maketitle

\begin{center}
\large{\textbf{0. Introduction}} \end{center} Let $k$ be a field and
$n$ be a natural number. By $G_n(k)=\mathrm{UT}(n, k)$ we denote the
group of all unipotent lower-triangular $n\times n$-matrices with
coefficients from $k$; this group is called a \emph{unitriangular}
group. By $\mathfrak{g}_n(k)=\mathfrak{ut}(n, k)$ we denote its Lie
algebra over $k$; this Lie algebra consists of all nilpotent
lower-triangular matrices with coefficients in $k$.

If $k=\mathbb{F}_q$ is a finite field, then $G_n(q)=G_n(k)$ is a
finite group; so, there are finitely many classes of equivalency of
irreducible complex representations of this group. A description of
all irreducible characters (or of some series of them) is a
classical problem of representation theory. The orbit method of
A.A.~Kirillov \cite{Kirillov1}, \cite{Kirillov2} allows to reduce
the similar problem of description of unitary irreducible
representations of Lie groups to the problem of classification of
coadjoint orbits; this method is also valid for $G_n(q)$ (see.
\cite{Kazhdan}), but complete classification of coadjoint orbits for
an arbitrary $n$ is unknown.

A description of orbits of the principal series (i.e., orbits of
maximum dimension) of Lie groups was presented in the pioneering
work on the orbit method \cite{Kirillov2}; it's also valid over a
finite field \cite{Kirillov3}. In C.~Andre's works (see
\cite{Andre1}, \cite{Andre2}) so-called \emph{basic characters} are
described (in particular, an exact formula for characters of the
principal series is found).

The problem of description of orbits, representations and characters
of sub-maximum di\-men\-si\-on is a natural generalisation of these
results; such orbits, representations and characters are called
\emph{subregular}. They play an important role in algebraic geometry
and $K$-theory (see, for example, \cite{Lusztig}). Subregular orbits
are described in \cite{IgnatevPanov}. The main goal of this paper is
to give an exact formula for the corresponding characters. This
formula shows that subregular characters (as characters of the
principal series) can be described in terms of coefficients of
minors of the characteristic matrix.

The paper is organized as follows. In section \ref{paragr_osn_ser}
we collect some basic facts about the group $G_n$ and discuss
Andre's formulae for characters of the principal series (Theorem
\ref{char_reg}). In section~\ref{paragr_sreg_formulir} we formulate
the Main Theorem including exact formulas for subregular characters
(Theorem \ref{theo_sreg}). More precisely, for an arbitrary
subregular character we find the elements such that the value of
this character on the conjugacy classes of these elements is
non-zero and compute this value. Section \ref{class_sopr} is devoted
to the description of these conjugacy classes (Theorem
\ref{theo_klass_sopr}). In section \ref{paragr_semidirect} we recall
some facts about semi-direct decomposition of the group $G_n$, wich
are needed for the proof of the main Theorem, which is given in
section \ref{paragr_proof}. Finally, section \ref{utochn_obob}
includes some remarks about discussed problems.
\par\bigskip

\sect{Characters of the principal series}\label{paragr_osn_ser}

Let $k=\mathbb{F}_q$ and $\mathrm{char}\,k=p$, i.e., $q=p^r$, where
$p$ is a prime number. Throughout the paper, we suppose $p\geqslant
n$. According to \cite{Kazhdan}, under these assumptions the orbit
method is valid: there is a one-to one correspondence between the
set of all irreducible complex characters of our group and the set
of all coadjoint orbits, i.e, $G_n(q)$-orbits in the dual space
$\mathfrak{g}_n^*(q)$. Moreover, to each orbit
$\Omega\subset\mathfrak{g}_n(q)$ the character
\begin{equation}\chi_{\Omega}(\mathrm{exp}\,a)
=q^{-\frac{1}{2}\dim\,\Omega}\cdot\sum_{f\in\Omega}
\theta(f(a)),\;a\in\mathfrak{g}_n(q),\label{tupaya_formula_char}\end{equation}
is assigned (here $\theta\colon\mathbb{F}_q\to\mathbb{C}$ is a
non-trivial character of the additive group of the field
$\mathbb{F}_q$). On the other hand, there are no exact formulas that
allow to find the value of the character of a given orbit on a given
element of the unitriangular group. \nota{There exists a form
$\langle A, B\rangle=\mathrm{tr}(AB)$ which is non-degenerate on
$\mathfrak{gl}_n(k)$. This allows to identify $\mathfrak{ut}^*(n,
q)$ with the space of all nilpotent upper-triangular matrices by the
formula $f(x)=\langle f, x\rangle=\sum_{i, j}\xi_{ji}x_{ij}$, where
$x=(x_{ij})\in\mathfrak{g}_n(q),\,f=(\xi_{ij})\in\mathfrak{g}_n^*(q)$.
Then the coadjoint action
$$K(g)\colon\mathfrak{g}^*_n\to\mathfrak{g}^*_n\colon(K(g)f)(x)=f(\mathrm{Ad}_
{g^{-1}}x),\,g\in G_n,\,f\in\mathfrak{g}_n^*,x\in\mathfrak{g}_n$$ is
given by $$K(g)\colon
x\mapsto(gxg^{-1})_{high},\,x\in\mathfrak{g}_n^*,g\in\mathfrak{g}_n,$$
where $(a)_{high}$ denotes the following matrix: its elements below
the main diagonal coincide with the corresponding elements of the
matrix $a$ and its elements on the main diagonal and below are equal
to zero.}

Representations of maximum dimension play an important role in this
theory. These re\-pre\-sen\-ta\-tions (and their orbits and
characters) form the so-called \emph{principal series}; orbits of
maximum dimension are called \emph{regular}. We'll use the following
\deno{Let $n\in\mathbb{N}$. Set $n_0=[n/2]$, $n_1=[(n-1)/2]$. Note
that $n=n_0+n_1+1$.} \deno{Let $g=(y_{ij})\in\mathrm{Mat}(n, k)$. By
$\Delta_{i_1,\ldots,i_k}^{j_1,\ldots,j_k}(g)$ we denote the minor of
the matrix $g$ with the rows $i_1,\ldots,i_k$ and the columns
$j_1,\ldots,j_k$ in the given order (for a given $1\lee k\lee
n,1\lee i_1,\ldots,i_k,j_1,\ldots,j_k\lee n$). In particular, for an
arbitrary $1\leqslant d\leqslant n_0$ we set $$\Delta_d(g)=
\Delta_{n-d+1,n-d+2,\ldots,n}^{1,2\ldots,d}(g).$$ We also denote
$\Delta^X(g)=\Delta_{\sigma(i_1,\ldots,i_k)}^{\tau(j_1,\ldots,j_k)}(g)$,
where $X$ is the set of the pairs $X=\{(i_1,j_1),\ldots,(i_k,j_k)\}$
and $\sigma, \tau$ are permutations such that
$\sigma(i_1)<\ldots<\sigma(i_k)$ and $\tau(j_1)<\ldots<\tau(j_k)$.}

A complete description of orbits of the principal series is
well-known \cite{Kirillov3}: \theo{Let $k=\mathbb{F}_q$ be a finite
field. An arbitrary regular orbit has the following defining
equations:
$$\Delta_d(^ta)=\beta_d,\;a\in\mathfrak{g}_n^*(q),\;1\leqslant d\leqslant n_0$$
(here $^ta$ is the transpose of the matrix $a$, $\beta_d$ are
arbitrary scalars from $\mathbb{F}_q$,
$\beta_1,\ldots,\beta_{n_0-1}\in\mathbb{F}_q^*$ and $\beta_{n_0}\neq
0$ for odd $n$).}\label{theo_reg}

\coro{For an arbitrary orbits of the principal series there exists
the unique \emph{canonical form of a regular orbit}, i.e., the
matrix of the form
$$f=\begin{pmatrix}
0&\ldots&0&\xi_{1, n}\\
0&\ldots&\xi_{2, n-1}&0\\
\vdots&\ddots&\vdots&\vdots\\
0&\ldots&0&0
\end{pmatrix},$$ where $\xi_{1, n}=\beta_1,\xi_{d, n-d+1}=
\dfrac{\beta_d}{\beta_{d-1}}, d=2,\ldots,n_0$, that consists in this
orbit.} \label{can_reg}

Moreover, dimension of a representation (and an orbit) of the
principal series is known \cite{Lehrer}: \theo{Let
$\mu(n)=(n-2)+(n-4)+\ldots$ and $T_{\Omega}$ be the representation
that corresponds to an orbit $\Omega$ of the principal series. Then
$\dim\Omega=2\mu(n),\;\dim T_{\Omega}=q^{\mu(n)}$.}

Note that there are exactly $q^{2\mu(n)}$ points on an orbit
$\Omega$ of dimension $2\mu(n)$ \cite{Kirillov4}.

In this section we recall the formula for characters of the
principal series, i.e., characters of the form $\chi=\chi_f$, where
$f\in\Omega_f\subset\g_n^*$ is the canonical form on a regular
orbit. More precisely, let $\chi$ be an irreducible character of the
group $G_n$. By $\mathrm{Supp}\mathop\chi$ denote its \emph{support}
(i.e., the set $\{g\in G_n\mid\chi(g)\neq0\}$); obviously, the
support is a union of certain conjugacy classes and the value of the
character on an arbitrary conjugacy class is constant. So, it's
enough to describe the support explicitly and compute the value of
the character on an arbitrary conjugacy class containing in the
support.

We'll use the following notations from Andre's paper \cite{Andre2}:
\deno{We denote by $\Phi(n)$ the set of all pairs $\{(i, j)\mid
1\lee j<i\lee n\}$ (we call them \emph{roots}, because to each
$(i,j)\in\Phi(n)$ the root vector $e_{ji}$, i.e., the matrix with 1
in the $(j,i)$-th entry and zeroes elsewhere, is assigned). Let
$D\subset\Phi(n)$ be a subset that contains at most one element from
each row and at most one element from each column; then this subset
is called \emph{basic}. If a basic subset consists from the roots of
the form $(n-j+1,j)$ then it's called \emph{regular}.}

\exam{Here we draw one of regular subsets: $D=\{(6, 1), ((4, 3))\}
\subset\Phi(6)$. The entries $(i, j)\in D$ are marked by the symbol
$\otimes$:
\begin{center}\small
$\mymatrix{
\pho & \pho&\pho& \pho& \pho& \pho\\
\Top{2pt}\Rt{2pt} \pho & \pho& \pho& \pho& \pho& \pho\\
\pho & \Top{2pt}\Rt{2pt} \pho& \pho& \pho& \pho& \pho\\
\pho & \pho& \Top{2pt}\Rt{2pt} \otimes& \pho& \pho& \pho\\
\pho & \pho& \pho& \Top{2pt}\Rt{2pt} \pho& \pho& \pho\\
\otimes & \pho& \pho& \pho& \Top{2pt}\Rt{2pt} \pho& \pho\\
}$
\end{center}}

\deno{For an arbitrary basis subset $D\subset\Phi(n)$ and arbitrary
map $\varphi\colon D\to\F_q^*$ consider the element of the group
$G_n(q)$ of the form $$x_D(\varphi)=1_n+\sum_{(n-j+1, j)\in
D}\varphi(n-j+1, j)e_{n-j+1, j},$$ where $1_n$ is the identity
$n\times n$-matrix (if $D=\varnothing$ then $x_D(\varphi)=1_n$). By
$\EuScript K_D(\varphi)$ denote the conjugacy class of this element
and by $\EuScript{K}_{\mathrm{reg}}$ denote the (disjoint) union of
$\EuScript K_D(\varphi)$ such that $D$ is a regular subset and
$\varphi\colon D\to\F_q^*$ is a map.} \defi{Let $D$ be a subset of
$\Phi(n)$. A root $(i,j)\in\Phi(n)$ is called $D$-\emph{regular}, if
$(i,k)\notin D$ and $(k,j)\notin D$ for any $i>k>j$. By $R(D)$
denote the set of all $D$-regular roots.}

\exam{Let $n=6$; $D=\{(3, 2), (6, 4)\}\subset\Phi(6)$ is a basic
subset. Here we mark the roots $(i, j)\notin R(D)$:
\begin{center}\small
$\mymatrix{
\pho & \pho&\pho& \pho& \pho& \pho\\
\Top{2pt}\Rt{2pt} \pho & \pho& \pho& \pho& \pho& \pho\\
\gray \pho & \Top{2pt}\Rt{2pt} \otimes& \pho& \pho& \pho& \pho\\
\pho & \gray \pho& \Top{2pt}\Rt{2pt} \pho& \pho& \pho& \pho\\
\gray \pho & \gray \pho& \gray \pho& \Top{2pt}\Rt{2pt} \otimes& \pho& \pho\\
\pho & \gray \pho& \pho& \gray \pho& \Top{2pt}\Rt{2pt} \pho& \pho\\
}$
\end{center}}

Now, set $\Phi_{\reg}=\{(i,j)\in\Phi(n)\mid i>n-j+1\}$,
$m_D=|R(D)\cap\Phi_{\reg}|$ for an arbitrary regular subset
$D\subset\Phi(n)$. Fix also any non-trivial character
$\theta\colon\F_q\to\mathbb{C}$ of the additive group of the ground
field.

\theop{Let $\Omega=\Omega_f\subset\g_n^*(q)$ be a regular orbit, let
$f=(\xi_{ij})$ be the canonical form of this orbit and $\chi=\chi_f$
be the corresponding character. Then
\begin{enumerate}
\item $\Supp{\chi}=\EuScript{K}_{\mathrm{reg}}$.
\item $\chi(g)=q^{m_D}\cdot\theta_f(e_D(\varphi))$ for any $g\in\EuScript
K_D(\varphi)\subset\EuScript{K}_{\mathrm{reg}}$, where
$e_D(\varphi)=x_D(\varphi)-1_n\in\g_n(q)$ and
$\theta_f\colon\g_n(q)\to\mathbb{C}$ is given by the formula
$$\theta_f(x)=\theta(f(x))=\prod_{(i,j)\in\Phi(n)}\theta(\xi_{ji}x_{ij}),
\quad x=(x_{ij})\in\g_n(q).$$ \end{enumerate}}{This is a special
case of \cite[Theorem 5.1]{Andre2}.}\label{char_reg}

In Andre's work \cite{Andre2} one can find explicit description of
$\EuScript K_D(\varphi)$ for any basic subset $D\subset\Phi(n)$.
More precisely, $g\in\EuScript K_D(\varphi)\subset G_n$ if and only
if
\begin{equation}
\Delta^{R_D(i, j)}(g)=\Delta^{R_D(i,
j)}(x_D(\varphi))\label{reg_uravn}
\end{equation}
for all $(i,j)\in R(D)$. Here $R_D(i, j)=\{(i, j)\}\cup\{(k, l)\in
D\mid l> j\text{ and }k< i\}$.

In particular, if $m=\max_{(i, j)\in D}j$ then
\begin{equation}
y_{ij}=0\text{ if }j>m\text{ or }i<n-m+1.\label{reg_zero}
\end{equation}
One can see that characters of the principal series and regular
orbits can be described in terms of minors of matrices from $G_n$;
by definition, this is also true for all Andre's basic characters.
In the next section we'll see that this is \emph{not} true for
subregular characters in general (see \cite{Andre2}).

\sect{Subregular characters
(statements)}\label{paragr_sreg_formulir} It follows  from Theorem
\ref{theo_reg} that a coadjoint orbit "in general position"\ has
maximum dimension $2\mu(n)$. More precisely, the set of points such
that their orbits are regular is a dense open subset of $\g_n^*(q)$
in Zariski topology. Denote this subset by $\mathcal O_{\reg}$, it's
given by the set of inequalities $\Delta_d\neq0$, $1\lee d\lee n_1$
(thats's why we use the term "regular orbits"). Recall that
$n=n_0+n_1+1$, where $n_0=[n/2]$, $n_1=[(n-1)/2]$.

On the other hand, for an arbitrary $1\lee d\lee n_1$ one can
consider the hyper-surface $\mathcal O_d\subset\g_n^*$ given by
$\Delta_d=0$. This hyper-surface splits into the union of certain
orbits.
\defi{An orbit $\Omega\in\mathfrak{g}_n^*$
(and the corresponding representation $T_{\Omega}$ and character
$\chi_{\Omega}$) is called \emph{subregular}, if
$\dim\Omega=2\mu(n)-2$ (resp., $\dim T_{\Omega}=q^{\mu(n)-1}$),
i.e., it has dimension $=$ dimension of a regular orbit $-\;2$ (in
the other words, this orbit has sub-maximum dimension, because any
orbit is even-dimensional \cite{Kirillov4}).}

Each subregular orbits contains in the unique $\mathcal O_d$ and has
maximum dimension among all orbits containing in $\mathcal O_d$.
Hence, we say that a subregular orbit is $d$-\emph{subregular} if it
contains in $\Omega\subset O_d$. In \cite{IgnatevPanov} defining
equations of subregular orbits are found, but we will not use them
in the sequel; we only list here elements of $\g^*_n$ such that
their orbits are subregular.
\defi{An element $f=(\xi_{ij})\in\g^*_n$ is called
\emph{a canonical form of a subregular orbit} (of the first, second
or third type resp.) if
\begin{enumerate}
\item There are a number $1\lee d<n_1$ and
$\beta_1,\ldots,\beta_{d-1},\beta',\beta'',\beta_{d+1}$, $\ldots$,
$\beta_{n_0-1}\in k^*$, $\beta_{n_0},\beta\in k$, where
$\beta_{n_0}\neq0$ for odd $n$, such that $\xi_{j,n-j+1}=\beta_j$
for any $1\lee j\lee d-1$ and $d+2\lee j\lee n_0$,
$\xi_{d,n-d}=\beta',\xi_{d+1,n-d+1}=\beta'',\xi_{n-d,n-d+1}=\beta$,
and $\xi_{ij}=0$ for all other roots.
\item $n$ is odd and there are $\beta_1,\ldots,\beta_{n_1-1}\in k^*,\beta',\beta''
\in k$ such that $\xi_{j,n-j+1}=\beta_j$ for any $1\lee j\lee
n_1-1$, $\xi_{n_1,n_0+1}=\beta',\xi_{n_1+1,n_0+2}=\beta''$,and
$\xi_{ij}=0$ for all other roots.
\item $n$ is even and there are $\beta_1,\ldots,\beta_{n_1-1},\beta\in k^*,
\beta',\beta''\in k$, such that $\xi_{j,n-j+1}=\beta_j$ for any
$1\lee j\lee n_1-1$ and either
$\xi_{n_1,n_0+1}=\beta,\xi_{n_1+1,n_0+2}=\beta',\xi_{n_1+2,n_0+2}=\beta''$
(and all other $\xi_{ij}=0$), or
$\xi_{n_1,n_0}=\beta',\xi_{n_1+1,n_0+2}=\beta$ (and all other
$\xi_{ij}=0$).
\end{enumerate}}\label{defi_sreg_orbits}

\exam{Here we draw some canonical forms of subregular orbits for
$n=8$, $d=1$ and $d=2$. The symbol $\otimes$ marks the roots $(i,
j)$ such that $\xi_{ji}\neq0$:
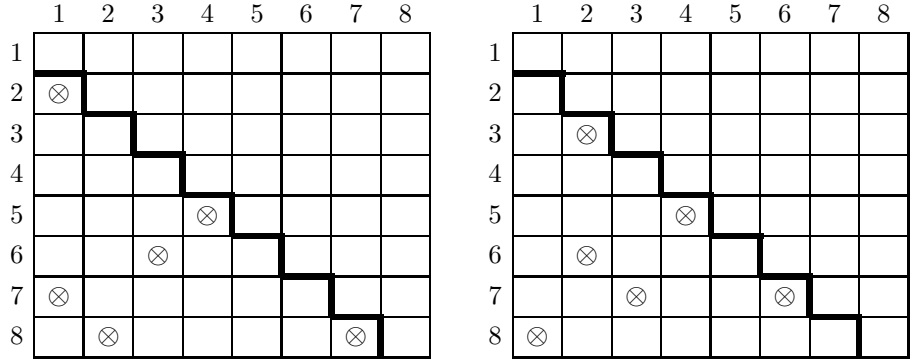
\begin{center}\small
$\mymatrix{
\pho & \pho&\pho& \pho& \pho& \pho &\pho &\pho\\
\Top{2pt}\Rt{2pt} \otimes & \pho& \pho& \pho& \pho& \pho &\pho &\pho\\
\pho & \Top{2pt}\Rt{2pt} \pho& \pho& \pho& \pho& \pho &\pho &\pho\\
\pho & \pho& \Top{2pt}\Rt{2pt} \pho& \pho& \pho& \pho &\pho &\pho\\
\pho & \pho& \pho& \Top{2pt}\Rt{2pt} \otimes& \pho& \pho &\pho &\pho\\
\pho & \pho& \otimes& \pho& \Top{2pt}\Rt{2pt} \pho& \pho &\pho &\pho\\
\otimes & \pho& \pho& \pho& \pho& \Top{2pt}\Rt{2pt} \pho &\pho &\pho\\
\pho & \otimes& \pho& \pho& \pho& \pho &\Top{2pt}\Rt{2pt} \otimes &\pho\\
}$\qquad$\mymatrix{
\pho & \pho&\pho& \pho& \pho& \pho &\pho &\pho\\
\Top{2pt}\Rt{2pt} \pho & \pho& \pho& \pho& \pho& \pho &\pho &\pho\\
\pho & \Top{2pt}\Rt{2pt} \otimes& \pho& \pho& \pho& \pho &\pho &\pho\\
\pho & \pho& \Top{2pt}\Rt{2pt} \pho& \pho& \pho& \pho &\pho &\pho\\
\pho & \pho& \pho& \Top{2pt}\Rt{2pt} \otimes& \pho& \pho &\pho &\pho\\
\pho & \otimes& \pho& \pho& \Top{2pt}\Rt{2pt} \pho& \pho &\pho &\pho\\
\pho & \pho& \otimes& \pho& \pho& \Top{2pt}\Rt{2pt} \otimes &\pho &\pho\\
\otimes & \pho& \pho& \pho& \pho& \pho &\Top{2pt}\Rt{2pt} \pho &\pho\\
}$
\end{center}}

\bigskip\theop{The orbit of a canonical form is subregular. More
over, for each subregular orbit there exists the unique canonical
form that contains in this orbit.}{ See \cite{IgnatevPanov} (there
the case $\mathrm{char}\,k=0$ is considered, but the proof is still
valid in our case $k=\F_q$ if $\chara k=p\gee n$).} This allows to
define $d$-subregular forms and subregular orbits of the first,
second or third type by the obvious way. Clearly, if $1\lee d<n_1$,
then $d$-subregular orbit is of the first type, and $n_1$-subregular
orbits are of the second type, if $n$ is odd, and the third one, if
$n$ is even.

So, the main goal is to find exact formulas for all characters of
the form $\chi_f$, where $f$ is a canonical form of a subregular
orbit. As for characters of the principal series, we'll describe the
support of a given subregular character and compute the value of
this characters on an arbitrary conjugacy class contained in the
support (we'll firstly consider the case $1\lee d<n_1$ and then (in
section \ref{utochn_obob}) the case of $n_1$-subregular orbits).

\deno{For an arbitrary $1\lee d< n_1$ by $D_0(d)$ and $D_1(d)$
denote one of the following sets resp.:
\begin{equation*}
\begin{split}
D_0(d):\quad&\varnothing,\{(n-d+1, d)\},\{(n-d,d)\},
\{(n-d+1,d+1)\},\\
&\{(n-d,d), (n-d+1,d+1)\},\\
D_1(d):\quad&\{(d+1,d), (n-d+1,n-d)\},\\
&\{(d+1,d), (n-d+1,n-d), (n-d,d), (n-d+1,d+1)\}.
\end{split}
\end{equation*}}

\defi{A subset $D\subset\Phi(n)$ is called
$d$-\emph{subregular} ($1\lee d<n_1$) if it has the form $D=D'\cup
D_i(d)$, where either $i=0$ or $i=1$ and $D'$ is a regular subset
that doesn't consist the roots $(n-d+1,d)$ and $(n-d+1, d+1)$. Note
that subregular subsets consisting $D_1(d)$ are not basic in
general.}{\label{defi_sreg_subset}} \nota{Denote $\Supp{f}=\{(i,
j\in\Phi(n))\mid f(e_{ij})\neq0\}$. One can see that
$D\subset\Supp{f}\cup(n-d+1,n-d)\cup(n-d+1,d)$, if $f$ is a
canonical form of a subregular orbit.}

\exam{Here we draw one of $1$-subregular subsets for $n=6$:

$D=\{(2, 1), (5, 1), (6, 2), (4, 3), (6, 5)\}\subset\Phi(6)$.
\begin{center}\small
$\mymatrix{
\pho & \pho&\pho& \pho& \pho& \pho\\
\Top{2pt}\Rt{2pt} \otimes & \pho& \pho& \pho& \pho& \pho\\
\pho & \Top{2pt}\Rt{2pt} \pho& \pho& \pho& \pho& \pho\\
\pho & \pho& \Top{2pt}\Rt{2pt} \otimes& \pho& \pho& \pho\\
\otimes & \pho& \pho& \Top{2pt}\Rt{2pt} \pho& \pho& \pho\\
\pho & \otimes& \pho& \pho& \Top{2pt}\Rt{2pt} \otimes& \pho\\
}$
\end{center}}

Define $R(D)$, $x_D(\varphi)$, $e_D(\varphi)$, $\EuScript
K_D(\varphi)$, $\theta$ and $\theta_f$, as for regular orbits. For a
subregular orbit $\Omega_f\subset\mathcal O_d$ denote
$$\EuScript K_f=\bigcup_{(D,\varphi)}\EuScript K_D(\varphi),$$ where
the union is over all $d$-subregular subsets $D$ and over maps
$\varphi\colon D\to\F_q^*$ such that
\begin{equation}\xi_{d,n-d}\cdot\varphi(d+1,d)=\xi_{d+1,n-d+1}
\cdot\varphi(n-d+1,n-d),\text{where
$f=(\xi_{ij})$}.\label{dop_usl}\end{equation}

For a $d$-subregular subset $D$ containing $D_1(d)$ we denote
\begin{equation*}
\Phi_d=\{(i,j)\in\Phi(n)\mid i>n-j+1, j\notin\{d, n-d\},i\notin
\{n-d+1, n-d\}\}.
\end{equation*}
Finally, for a $d$-subregular subset $D$ let
\begin{equation*}
m_D=\begin{cases}|R(D)\cap\Phi_{\reg}|-1,&\text{if }D\supset
D_0(d),\\
|R(D)\cap\Phi_d|+n-2d-1,&\text{if }D\supset D_1(d).
\end{cases}
\end{equation*}

\theosp{Let $1\lee d<n_1$, let $\Omega_f\subset\mathcal
O_d\subset\g_n^*(q)$ be a subregular orbit, and $\chi=\chi_f$ be the
corresponding character. Then
\begin{enumerate}
\item$\Supp{\chi}=\EuScript K_f$.
\item $\chi(g)=q^{m_D}\cdot\theta_f(e_D(\varphi))$
for an arbitrary $g\in\EuScript K_D(\varphi)\subset\EuScript K_f$.
\end{enumerate}}{In Lemma \ref{lemm_proofs_inv}
we present certain equations the element $x_D(\varphi)$ satisfies,
and show that the ideal $J$ of $k[G_n]$ generated by these equations
is invariant under the adjoint action. In Lemma
\ref{lemm_proofs_prim} we prove that $J$ is a prime ideal. In Lemma
\ref{lemm_centr} we find the stabilizer $\mathcal
C=\mathrm{Stab}_{\,G_n}(x_D(\varphi))$ and prove that
$$\dim\EuScript K_D(\varphi)=\codim{\mathcal C}=\dim V(J)$$
(here $V(J)=\{g\in G_n\mid F(g)=0\text{ for all }F\in J\}$, as
usual). This shows that $\EuScript K_D(\varphi)=V(J)$, i.e., our
equations are exactly the defining equations for the conjugacy class
of the element $x_D(\varphi)$ (Theorem \ref{theo_klass_sopr}).

On the other hand, in Lemmas \ref{lemm_O_1_neq} ---
\ref{lemm_for_all_d} we prove that if $\chi(g)\neq0$, then $g$
satisfies these equations (and the additional equation
(\ref{dop_usl})), and compute the value $\chi(g)$ in this case. This
concludes the proof.}\label{theo_sreg}

\nota{If $D\supset D_1(d)$ is $d$-subregular, $\vfi, \w\vfi$ are
maps from $D$ to $\F_q^*$, and
\begin{equation*}\begin{split}&\vfi(d+1, d)\cdot\vfi(n-d+1, d+1)+
\vfi(n-d, d)\cdot\vfi(n-d+1,n-d)=\\
&=\w\vfi(d+1, d)\cdot\w\vfi(n-d+1, d+1)+\w\vfi(n-d,
d)\cdot\w\vfi(n-d+1, n-d),
\end{split}\end{equation*}
then $\K_D(\vfi)=\K_D(\w\vfi)$ (see section \ref{class_sopr}). But
one can see (using the results of section \ref{class_sopr}) that in
this case $\theta_f(e_D(\vfi))=\theta_f(e_D(\w\vfi))$. Since $m_D$
is independent of $\vfi$, the character value in the formulation of
the Theorem is well-defined.}

It's easy to prove a variant of this Theorem for the case of
$n_1$-subregular orbits by some modifications of definitions and
formulations (see section \ref{utochn_obob}).

\sect{Conjugacy classes $\K_D(\varphi)$}\label{class_sopr} Here we
present an explicit description of the conjugacy class
$\K_D(\varphi)$ of an element $x_D(\varphi)$. Fix an arbitrary
$1\lee d<n_1$ and a subregular subset $D\supset D_1(d)$ (in this
section and in the two next sections we consider the case
$D(d)=D_1(d)$; the case $D(d)=D_0(d)$ is similar, see
section~\ref{utochn_obob}). \deno{It's convenient to split a
$d$-subregular subset into the union
$$D=D^-\sqcup D_1(d)\sqcup D^+,$$ where $D^-=\{(i,j)\in D\mid j<d\}$,
$D^+=\{(i,j)\in D\mid d<j<n-d\}$.

Recall that $D'=D\setminus D_1(d)=D^-\sqcup D^+$ is a regular subset
(see Definition \ref{defi_sreg_subset}). Let $D''=D\setminus\{(n-d,
d), (n-d+1, d+1)\})$ (it's a basic subset of $\Phi(n)$, and
$D'\subset D''\subset D$).}

It's impossible to describe $\EuScript K_D(\varphi)$ in terms of
minors of matrices $g=(y_{ij})\in G_n$. Let $(i, j)\in D^+$ and
$m=\max_{(i, j)\in D^+}j$. Consider the following polynomials:
\begin{equation*}
\alpha_{ij}=\sum_{l=n-m+1}^{n-d}y_{n-d+1,l}y_{l,j},\quad \beta_{ij}=
\sum_{l=d+1}^{m}y_{i,l}y_{l,d},\quad
\gamma=\sum_{l=d+1}^{n-d}y_{n-d+1,l}y_{l,d}
\end{equation*}
(more conceptual description of these polynomials is given in
section \ref{utochn_obob}). For simplicity, we'll write
$\Delta_{ij}$ instead of $\Delta^{R_{D''}(i,j)}(g)$, and $y_{\beta},
y_{\alpha}$ instead of $y_{d+1,d},y_{n-d+1,n-d}$ resp.
($R_{D''}(i,j)$ are defined similar to $R_D(i, j)$ in
(\ref{reg_uravn})).

Let us now start with proving Theorem \ref{theo_sreg}. Let $c_0\in
k$, $c_{\alpha}, c_{\beta},c_{ij}\in k^*$, where $(i,j)\in D'$, are
arbitrary scalars. Consider the ideal $J$ of $k[G_n]$ generated by
the elements
\begin{equation}\begin{split}
&y_{\alpha}-c_{\alpha}, y_{\beta}-c_{\beta},\gamma-c_0,\\
&\alpha_{ij},\beta_{ij},\quad(i,j)\in D^+,\\
&\Delta_{ij}-c_{ij},\quad(i,j)\in D',\\
&\Delta_{ij},\quad(i,j)\in R(D'')\setminus D'.
\end{split}\label{sreg_uravn}\end{equation}
\lemmp{$J$ is $G_n$-invariant, i.e., if $g\in V(J)$, then
$xgx^{-1}\in V(J)$ for all $x\in G_n$.}{Since the set $D''$ is
basic, all $\Delta_{ij}$ invariant \cite[Lemma 2.1]{Andre2}.

Let $g=(y_{ij})\in V(J)$. Denote $x_{rs}(\lambda)=1+\lambda e_{rs}$.
Since every element $x\in G_n$ can be written as a product
$x=x_{r_1s_1}(\lambda_1)\ldots x_{r_ms_m}(\lambda_m)$ for a certain
$m$ and $(r_i,s_i)\in\Phi(n)$, $\lambda_i\in k$, it's enough to
prove that if $x=x_{rs}(\lambda)$, then $xgx^{-1}\in V(J)$. But if
$x$ has this form, then
\begin{equation*} (xgx^{-1})_{ij}=\begin{cases} y_{ij},&\text{if
}i\neq r\text{ and
}j\neq s,\\
y_{ij},&\text{if }i=r\text{ and }j\gee s\text{, or }j=s\text{ and }i\lee r,\\
y_{rj}+\lambda y_{sj},&\text{if }i=r,j<s,\\
y_{is}-\lambda y_{ir},&\text{if }j=s,i>r.
\end{cases}
\end{equation*}
Hence, $y_{\alpha},y_{\beta}$ are invariant obviously, and the proof
of invariance of other elements is by direct enumeration of possible
values of $r$ and $s$.

For example, consider the polynomial $\alpha_{ij}$. If $s\lee j$ or
$r> n-d+1$, then all coordinate functions involved in this
polynomial are invariant themselves; this is also true, if
$r=n-d+1$, $s\lee n-m+1$. If $r=n-d+1,s> n-m+1$, then
\begin{equation*}
\begin{split}
\alpha_{ij}(xgx^{-1})&=\sum_{l=n-m+1}^{s-1}(y_{n-d+1,l}+\lambda
y_{s,l})y_{l,j}+
\sum_{l=s}^{n-d}y_{n-d+1,l}y_{l,j}=\\
&=\alpha_{ij}(g)+\lambda\cdot\sum_{l=n-m+1}^{s-1}
y_{s,l}y_{l,j}=\alpha_{ij}(g),
\end{split}
\end{equation*} because
$y_{s,l}=0$ for all $l\gee n-m+1\gee n-n_0+1>m$ (it's a particular
case of (\ref{reg_zero})).

If $r<n-d+1,s= j$, then $\alpha_{ij}$ is invariant for the same
reasons. Finally, if $r<n-d+1,s> j$, then
\begin{equation*}\begin{split} \alpha_{ij}(xgx^{-1})&= \sum_{l\neq
r,s}y_{n-d+1,l}y_{l,j}+(y_{r,j}+\lambda
y_{s,j})y_{n-d+1,r}+y_{s, j}(y_{n-d+1,s}-\lambda y_{n-d+1,r})=\\
&=\alpha_{ij}(g)+\lambda\cdot(y_{s,j}y_{n-d+1,r}-y_{s,
j}y_{n-d+1,r})=\alpha_{ij}(g).
\end{split}\end{equation*}
Invariance of $\beta_{ij}$ and $\gamma$ can be proved similarly.}
\label{lemm_proofs_inv}

This means that $V(J)$ is a union of cojugacy classes. It's easy to
see that $x_D(\varphi)\in V(J)$, if
\begin{equation}
\begin{split}
&c_{ij}=\Delta_{ij}(x_D(\varphi)),\\
&c_0=\gamma(x_D(\varphi)),\\
&c_{\alpha}=\varphi(n-d+1,n-d),\\
&c_{\beta}=\varphi(d+1,d).\end{split}\label{const_x_D}
\end{equation} Thus, $\K_D(\varphi)\subset V(J)$ for this values
of constants.

For any root $\xi=(i,j)\in\Phi(n)$, we define its \emph{level} as
the number $u(\xi)=i-j$. Then the formula
\begin{equation} \xi=(i_1,j_1)<\eta=(i_2,j_2)\Leftrightarrow
\text{ either }u(\xi)<u(\eta)\text{ or }u(\xi)=u(\eta),j_1<j_2,
\end{equation}
defines a complete order on the set of all roots. For each root
$\xi=(i,j)\in\Phi(n)$ let $I_{\xi}$ be the ideal in $k[G_n]$
generated by all $y_{\eta},\eta<\xi$, and $\xi_0=(i_0,j_0)$ be the
maximal root in $D''$, which is less than $\xi$ (if exists).

\newpage
\lemmp{The ideal $J$ is a prime ideal of $\bar k[y_{ij}]$ (here
$\bar k$ is the algebraic closure of $k$).}{Consider the following
transformation of coordinates:
\begin{equation}
\begin{split}
&\w y_{n-d+1,n-d}=y_{\alpha}-c_{\alpha},\quad \w y_{d+1,d}=
\w y_{\beta}-c_{\beta},\quad \w y_{n-d+1,d+1}=\gamma-c_0,\\
&\w y_{n-d+1,i}=\alpha_{ij},\quad
\w y_{i,d+1}=\beta_{ij},\quad(i,j)\in D^+,\\
&\w y_{ij}=\Delta_{ij}-c_{ij},\quad(i,j)\in D',\\
&\w y_{ij}=\Delta_{ij},\quad(i,j)\in R(D'')\setminus D'.
\end{split}\label{zamena}
\end{equation}
(it's easy to see that $J=\langle \w y_\xi\rangle_{\xi\in B}$, where
$B\subset\Phi(n)$ denotes the set of all roots $(i,j)$ from the
left-side hand of (\ref{zamena})). Note that for any $\xi\in B$ we
have
\begin{equation}\w y_{\xi}\equiv y_{\xi}^0\cdot
y_{\xi}+a_{\xi}\pmod{I_{\xi}},\label{treug}\end{equation} where
$y_{\xi}^0$ is an invertible element of $\bar k[y_{ij}]/J$, and
$a_{\xi}\in k$ is a certain scalar. Indeed, this is evident for $\w
y_{n-d+1,n-d}$, $\w y_{d+1,d}$. For other roots we have:
\begin{equation*}
\begin{split}
&\w y_{n-d+1,d+1}=\gamma-c_0=y_{n-d+1,d+1}y_{d+1,d}+\ldots,\\
&\w y_{n-d+1,i}=\alpha_{ij}= y_{n-d+1,i}y_{i,j}+\ldots,\quad(i,j)\in
D^+,\\
&\w y_{i,d+1}=\beta_{ij}=y_{i,d+1}y_{d+1,d}+\ldots,\quad(i,j)\in D^+,\\
&\w y_{ij}=\Delta_{ij}-c_{ij}=y_{ij}\cdot\Delta_{i_0,j_0}+\ldots,\quad(i,j)\in D',\\
&\w
y_{ij}=\Delta_{ij}=y_{ij}\cdot\Delta_{i_0,j_0}+\ldots,\quad(i,j)\in
R(D'')\setminus D',
\end{split}
\end{equation*}
where for any $\w y_{\xi}$, $\xi\in B$, dots denote elements, equal
to zero modulo $I_{\xi}$, and scalars (we assume that
$\Delta_{i_0,j_0}=1$, if the root $\xi_0$ does not exist for a given
$\xi\in\Phi(n)$). But one can easily obtain the following equatities
modulo $J$:
\begin{equation*}\begin{split} &y_{d+1,d}\equiv
c_{\beta}\neq0,\quad y_{n-m+1,m}\equiv c_{n-m+1,m}\neq0,\\
& y_{i,j}\equiv c_{ij}/c_{i_0,j_0}\neq0,\quad(i,j)\in D^+,j<m,\\
&\Delta_{i_0,j_0}\equiv c_{i_0,j_0}\neq0,\quad(i_0,j_0)\in D''
\end{split}
\end{equation*}
(recall that $m=\max_{(i,j)\in D^+}j$). This concludes the proof of
(\ref{treug}).

Hence, in $\bar k[y_{ij}]/J$ all $y_{\xi}$ are polynomial in $\w
y_{\xi}$ (for each $\xi\in B$). Consequently,
$$\bar k[y_{ij}]/J=\left.\bar k[y_{\xi}]_{\xi\in\Phi(n)}
\right/\langle\w y_{\xi}\rangle_{\xi\in B} \cong\left.\bar
k[\{y_{\xi}\}_{\xi\notin B}\cup\{\w y_{\xi}\}_{\xi\in
B}]\right./\langle\w y_{\xi}\rangle_{\xi\in B}\cong \bar k[\w
y_{\xi}]_{\xi\notin B}.$$ In particular, $\bar k[y_{ij}]/J$ is a
domain, so $J$ is a prime ideal.}\label{lemm_proofs_prim}

This Lemma shows that $V(J)$ is an irreducible subvariety in
$G_n(\bar k)$, because $\bar k[G_n]\cong\bar k[y_{ij}]$. \lemmp{Let
$\mathcal C=\mathrm{Stab}_{\,G_n}(x_D(\varphi))=\{g\in G_n\mid
gx_D(\varphi)=x_D(\varphi)g\}$ be the stabilizer (the centralizer)
of $x_D(\varphi)$. Then $\dim \mathcal C=\codim V(J)$, where $J$ is
generated by the elements (\ref{sreg_uravn}) with the scalars given
by (\ref{const_x_D}).}{Firstly, we'll present the defining equations
of the centralizer $\mathcal C$. For any $\xi=(i,j)\in\Phi(n)$, let
$\Phi_{\xi}$ be the union of all roots form the $i$-th column and
the $j$-th row. We put
\begin{equation*}\begin{split}
&\Phi_{\alpha}=\{(n-d+1,i)\mid(i,j)\in
D^+\},\quad\w\Phi_{\alpha}=\{(n-d+1,j)\mid j<d\},\\
&\Phi_{\beta}=\{(i,d+1)\mid(i,j)\in
D^+\},\quad\w\Phi_{\beta}=\{(i,d)\mid i>n-d+1\},\\
&\Phi_{\gamma}=\{(n-d,d)\},\quad\Phi_{\delta}=(\cup_{\xi\in D}\Phi_{\xi})\cap (R(D)\setminus D),\\
&A=\Phi_{\alpha}\sqcup\w\Phi_{\alpha}\sqcup
\Phi_{\beta}\sqcup\w\Phi_{\beta}\sqcup\Phi_{\gamma}\sqcup\Phi_{\delta}\subset\Phi(n)
\end{split}\end{equation*}
(these sets are really disjoint). It's easy to check that $\mathcal
C$ is given by the following equations:
\begin{equation}\begin{split}
&(\alpha)\quad y_{n-d+1,i}a_{ij}=a_{n-d+1,n-d}y_{n-d,j},\quad(i,j)\in D^+,\\
&(\beta)\quad y_{i,d+1}a_{d+1,d}=a_{ij}y_{jd},\quad(i,j)\in D^+,\\
&(\w\alpha)\quad a_{n-d+1,n-d}y_{n-d,j}+a_{n-d+1,d+1}y_{d+1,j}=0,
\quad j<d,\\
&(\w\beta)\quad y_{i,d+1}a_{d+1,d}+y_{i,n-d}a_{n-d,d}=0,\quad i>n-d+1,\\
&(\gamma)\quad y_{n-d+1,n-d}a_{n-d,d}+y_{n-d+1,d+1}a_{d+1,d}=
a_{n-d+1,n-d}y_{n-d,d}+a_{n-d+1,d+1}y_{d+1,d},\\
&(\delta)\quad y_{ij}=0,\quad(i,j)\in\Phi_{\delta},\\
\end{split}\label{centralizer}\end{equation}
(here we write $a_{ij}=\varphi(i,j)\in k^*$ for simplicity).

So, $\mathcal C$ is defined by the equations (\ref{centralizer}),
which are labeled  by roots from $A$, and $V(J)$ is defined by the
equations (\ref{sreg_uravn}), which are labeled by roots from $B$.
Consider the maps $\sigma_A\colon A\to\Phi(n)$ and $\sigma_B\colon
B\to\Phi(n)$, given by the formulas

\begin{equation*}
\begin{split}
&\sigma_A(\xi)=\begin{cases}(n-i+1,j),&
\xi=(i,j),(n-i+1,i)\in D'',j\neq d+1,i<n-j+1,\\
(i,n-j+1),&
\xi=(i,j),(n-j+1,j)\in D'',i\neq n-d,i>n-j+1,\\
(i,d),&
\xi=(i,d+1)\in\Phi_{\delta},i\neq n-d,\\
(n-d+1,j),&
\xi=(n-d,j)\in\Phi_{\delta},j\neq d+1,\\
(n-d+1,d),&
\xi=(n-d,d+1)\in\Phi_{\delta},\\
\xi,&\text{elsewhere},
\end{cases}\\
&\\
&\sigma_B(\xi)=\begin{cases} (n-d,j),&
\xi=(n-d+1,i),(i,j)\in D^+,\\
\xi,&\text{elsewhere}.
\end{cases}
\end{split}
\end{equation*}
They define the map $\sigma\colon A\sqcup B\to\Phi(n)$ (here
$A\sqcup B$ is the disjoint union of $A$ and $B$) by the rule
$\sigma\mathop{\mid}_A=\sigma_A$, $\sigma\mathop{\mid}_B=\sigma_B$;
it's easy to see that this map is a bijection. Since
$\dim\K_D(\varphi)=\codim\mathcal C=|A|$, $\codim V(J)=|B|$, this
concludes the proof.}\label{lemm_centr}

\theop{The defining ideal $J$ of the conjugacy call
$\K_D(\varphi)\subset G_n(k)$ of an element $x_D(\varphi)$ is
generated by the elements (\ref{sreg_uravn}) with the scalars given
by (\ref{const_x_D}).}{Since conjugacy classes of $G_n(\bar k)$ are
Zariski-closed \cite[Proposition 2. 5]{Steinberg}, previous Lemmas
follow that $\K_D(\varphi)=V(J)$ over $\bar k$. Hence, the sets of
their $k$-points also coincide.}\label{theo_klass_sopr}

\sect{A semi-direct decomposition of $G_n$}\label{paragr_semidirect}
The next goal is to describe the support of a subregular character
$\chi=\chi_f$ explicitly (we'll see that the support coincides with
$\K_f$) and to compute the value of this character on a counjugacy
class contained in the support. We'll use induction on dimension of
the group and the Mackey's method of semi-direct decomposition (see,
f.e., \cite{Lehrer}). In this section we collect some basic facts
which are needed for the sequel.

Let $G$ be a finite group, $A, B$ be its subgroups and $G=A\rtimes
B$ be their semi-direct product (i.e., $G=AB$ and $A\triangleleft
G$).
\defi{Let $G=A\rtimes B$ be a finite group and $A$ be abelian. For a given
irreducible character $\psi$ of the group $A$, the subset
$B^{\psi}=\{b\in B\mid\psi\circ\tau_b=\psi\}$ of the group $B$ is
said to be the \emph{centralizer} of this character (here
$\tau_b\colon A\to A\colon a\mapsto bab^{-1}$).}

The following Theorem satisfies \cite{Lehrer}: \theo{Let $G=A\rtimes
B$ be a finite group and $A$ be abelian. Then every irreducible
representation $\tau$ of the group $G$ has the form
$\tau=\Ind{A\rtimes B^{\chi}}{G}{\psi\otimes\widetilde\tau}$, where
$\psi$ is a certain irreducible character of the group $A$, and
$\w\tau$ is a certain irreducible representation of the centralizer
$B^{\psi}$. Hence, every irreducible character $\chi$ of the group
$G$ has the form $\chi=\Ind{A\rtimes
B^{\psi}}{G}{\psi\widetilde\chi}$, where $\widetilde\chi$, $\psi$
are certain irreducible characters of the groups $A$, $B_{\psi}$
resp. On the other hand, any character $\Ind{A\rtimes
B^{\psi}}{G}{\psi\widetilde\chi}$ is an irreducible character of
$G$.}\label{semi_direct}

We denote \begin{equation*} \begin{split}&P_n=\{g=(y_{ij})\in
G_n\mid
y_{ij}=0\text{ for }j\neq1\},\\
&G_{n-1}\cong\{g=(y_{ij})\in G_n\mid y_{ij}=0\text{ for
}j=1\}\hookrightarrow G_n.\\
\end{split}
\end{equation*}

Then $G_n=P_n\rtimes G_{n-1}$, moreover, the group
$P_n\cong\F_q^{n-1}$ is abelian; hence, all conditions of Theorem
\ref{semi_direct} are satisfied.

We fix a non-trivial additive character
$\theta\colon\F_q\to\mathbb{C}$. Any irreducible character of the
group $P_n$ has the form
\begin{equation*}
p=(p_{ij})\in P_n\mapsto\theta(s_2 p_{21})\cdot\ldots\cdot\theta(s_n
p_{n1}),
\end{equation*}
where $s_i\in\F_q$, $2\lee i\lee n$, are arbitrary scalars. We'll
consider the case $\theta_n(p)=\theta(s_np_{n1})$ and the case
$\theta_{n-1}(p)=\theta(s_{n-1}p_{n-1,1})$. One can see that their
centralizers in the subgroup $G_{n-1}$ are
\begin{equation*}
\begin{split}
&G_{n-1}^{\theta_n}=\{g=(y_{ij})\in G_{n-1}\subset G_n\mid
y_{nj}=0\text{ for }1\lee j\lee n-1\},\\
&G_{n-1}^{\theta_{n-1}}=\{g=(y_{ij})\in G_{n-1}\subset G_n\mid
y_{n-1,j}=0\text{ for }1\lee j\lee n-2\}.
\end{split}
\end{equation*}
Note that $G_{n-1}^{\theta_n}\cong G_{n-2}$, and
$G_{n-1}^{\theta_{n-1}}\cong\w G_{n-2}=G_{n-2}\times\F_q$ (we'll
consider only these embeddings of these subgroups into $G_n$).

For a linear function $f\in\g_n^*$, let $\pi(f)\in\g_{n-2}^*$ and
$\w\pi(f)$ denote its restrictions to $\g_{n-2}$ and $\w\g_{n-2}$
respectively (here $\w\g_{n-2}=\mathop{\mathrm{Lie}}\w G_{n-2}$ and
embeddings of subalgebras into $\g_n$ correspond to embeddings of
subgroups into $G_n$). Due to the Mackey's method, if $1\lee d<n_1$
and $f=(\xi_{ij})\in\mathcal O_d$ is the canonical form of a
$d$-subregular orbit, then
\begin{equation}
\begin{split}
&T_f=\Ind{P_n\rtimes G_{n-2}}{G_n}{\theta_n\otimes
T_{\pi(f)}},\text{ where }s_n=\xi_{1n},\quad\text{
if }d>1,\\
&T_f=\Ind{P_n\rtimes\w G_{n-2}}{G_n}{\theta_{n-1}\otimes
T_{\w\pi(f)}},\text{ where }s_{n-1}=\xi_{1,n-1},\quad\text{
if }d=1,\\
\end{split}\label{ind_reps}
\end{equation}
where $T_f$ (resp. $T_{\pi(f)}$ and $T_{\w\pi(f)}$) denotes the
representation of the group $G_n$ (resp. $G_{n-2}$ and $\w
G_{n-2}$), corresponding to the orbit $\Omega_f\subset\g_n^*$ (resp.
$\Omega_{\pi(f)}\subset\g_{n-2}^*$ and
$\Omega_{\w\pi(f)}\subset\w\g_{n-2}^*$). All terms in the last
equality are well-defined, because the orbit method is valid for the
group $\w G_{n-2}$. So, the problem can be reduced to the study of
representations of dimension less than the first one has, and we can
use an inductive argument.

Finally, certain coset decompositions of the group $G_n$ are needed
for construction of induced representations. It's easy to check that
\begin{equation*}
\begin{split}
& H_n=\{h=(t_{ij})\in G_{n-1}\mid t_{ij}=0\text{ for }i\neq n\},\\
& \w H_n=\{h=(t_{ij})\in G_{n-1}\mid t_{ij}=0\text{ for }i\neq
n-1\}.
\end{split}
\end{equation*} are complete systems of representatives of
$G_n/(P_n\rtimes G_{n-2})$ and $G_n/(P_n\rtimes\w G_{n-2})$ resp.
\nota{Since $\w G_{n-2}\cong G_{n-2}\times\F_q$, we have a complete
description of the representation $T_{\w\pi(f)}$. Precisely, we note
that $G_{n-2}$ is isomorphic to the subgroup $G_{n-2}\times0$ of $\w
G_{n-2}$ (this induces the embedding of Lie algebras
$\g_{n-2}\hookrightarrow\w\g_{n-2}$); let $\psi$ be the projection
$\w G_{n-2}\to G_{n-2}=G_{n-2}\times0$ and $g=(y_{ij})\in\w
G_{n-2}$. Then the character of the representation $T_{\w\pi(f)}$
has the form
$g\mapsto\chi(\psi(g))\cdot\theta(\xi_{n-1,n}y_{n,n-1})$, where
$\chi$ is the character of the principal series of $G_{n-2}$,
corresponding to the orbit
$(\w\pi(f))\mathbin{\mid}_{\g_{n-2}}$.}\label{nota_w_pi} \nota{There
is another ("symmetric") decomposition $G_n=P_n'\rtimes G_{n-1}'$,
where
\begin{equation*}
\begin{split}
&P_n'=\{g=(y_{ij})\in G_n\mid y_{ij}=0\text{ for }j\neq1\},\\
&G_{n-1}'=\{g=(y_{ij})\in G_n\mid y_{ij}=0\text{ for }j=1\}.
\end{split}
\end{equation*}
Reflecting other subgroups and subsets in the anti-diagonal, we get
$\w G_{n-2}'$, $H_n'$ and $\w H_n'$ ($G_{n-2}$, embedded into $G_n$
as above, is invariant under this reflection, i.e,
$G_{n-2}'=G_{n-2}$). Then the irreducible character
$\theta_{n-1}'\colon P_n'\to\mathbb{C}$ has the form $P_n'\ni
p=(p_{ij})\mapsto\theta(\xi_{2,n}p_{n,2}+\xi_{n-1,n}p_{n,n-1})$.}
\label{nota_symm_rtimes} Now we are able to conclude the proof of
Theorem \ref{theo_sreg}.

\sect{The proof of the Main Theorem}\label{paragr_proof} Let $1\lee
d<n_1$, $f\in\mathcal O_d$ be a canonical form of a $d$-subregular
orbit and $\chi=\chi_f$ be the corresponding character (as above).
We'll prove in this section that its support coincides with $\K_f$
and compute its value on an arbitrary conjugacy class
$\K_D(\varphi)\subset\K_f$. We've proved in section~\ref{class_sopr}
that the subvariety $\K_D(\varphi)$ of $G_n$ is defined by the
equations (\ref{sreg_uravn}) with the scalars given by
(\ref{const_x_D}); in the other words, $V(J)$ is the defining ideal
of $\K_D(\varphi)$ in $k[G_n]$.

The following proof is by induction on $n$; the base can be checked
directly (f.e., using (\ref{tupaya_formula_char})). We'll assume
that if $g=(y_{ij})\in G_n$, then $y_{n-d+1,n-d}\neq0$ (this means
that $D\supset D_1(d)$). Firstly, we have to consider the case
$d=1$. For convenience, we denote
\begin{equation*}
\begin{split}
&\Phi_{i_0}=\{(i,j)\in\Phi(n)\mid i=i_0\},\quad
\w\Phi_{i_0}=\cup_{i\gee i_0}\Phi_i,\\
&\Phi^{j_0}=\{(i,j)\in\Phi(n)\mid j=j_0\},\quad
\w\Phi^{j_0}=\cup_{j\lee j_0}\Phi^j.\\
\end{split}
\end{equation*}
\lemmp{Let $f\in\mathcal O_1$ be a canonical form of a
$1$-subregular orbit, and $g\in G_n$. If $\chi(g)\neq0$, then
$g\in\K_D(\varphi)$ for a certain $x_D(\varphi)$ satisfying
(\ref{dop_usl}).}{Recall that
$$T_f=\Ind{P_n\rtimes\w G_{n-2}}{G_n}{\theta_{n-1}\otimes
T_{\w\pi(f)}}$$ (see (\ref{ind_reps})). Since $G_n=P_n\rtimes
G_{n-1}$, an arbitrary element $g\in G_n$ can be uniquely
represented as $g=pg'$, $p\in P_n$, $g'\in G_{n-1}$ (in fact, $p$
and $g'$ are given by replacing the corresponding entries of $g$ by
zeroes). Hence, for a given $g\in G_n$, the element $p(g)\in P_n$ is
well-defined. By $\w\pi(g)$ we denote the element of $\w G_{n-2}$,
which is given by replacing the corresponding entries of $g$ by
zeroes. We have
\begin{equation}
\chi(g)=\chi_f(g)=\Ind{P_n\rtimes\w
G_{n-2}}{G_n}{\theta_{n-1}(p(g))\cdot\chi_{\w\pi(f)}(\w\pi(g))}=\sum_{h\in\w
H_n}\theta_{n-1}(p(h^{-1}gh))\cdot\chi_{\w\pi(f)}(\w\pi(h^{-1}gh))
\label{first_sreg_char_formula}
\end{equation}
(here the summation is over all $h\in\w H_n$ such that $h^{-1}gh\in
P_n\rtimes\w G_{n-2}$).

Note that for any $h=(t_{ij})\in\w H_n$,
\begin{equation}
\begin{split}&(h^{-1}gh)_{n,j}=y_{n,j}+y_{n,n-1}t_{n-1,j},
\quad\text{if $2\lee j\lee n-2$},\\
&(h^{-1}gh)_{n-1,j}=y_{n-1,j}-\sum_{i=j+1}^{n-2}t_{n-1,i}y_{i,1},\quad\text{if
$1\lee j\lee n-3$},\\
&(h^{-1}gh)_{ij}=y_{ij}\quad\text{for all other $1\lee j<i\lee n$}.
\end{split}{\label{hgh}}\end{equation}
If $\chi_{\w\pi(f)}(\w\pi(g))\neq0$, then it follows from the
equations of regular orbits (Theorem \ref{theo_reg}), formulas for
characters of the principal series (Theorem \ref{char_reg}) and
remark \ref{nota_w_pi} that
\begin{equation}
\chi_{\w\pi(f)}(\w\pi(g))=q^s\cdot\chi_{\pi^2(f)}(\pi^2(g))\cdot\theta(\xi_{2,n}
\cdot(-1)^{|X_1|-1}\cdot\frac{\Delta^{X_1}(g)}{\Delta^{X_2}(g)})
\cdot\theta(\xi_{n-1,n}y_{n,n-1}). \label{second_sreg_char_formula}
\end{equation}
Here $X_1=D\setminus(\Phi^1\cup\Phi^{n-1})$,
$X_2=X_1\setminus\Phi_n$,
$s=|(R(D^+)\cap\Phi_n)\setminus(\w\Phi^2\cup\Phi^{n-1})|$ and
$\pi^2$ denotes the maps $G_n\to G_{n-4}$ and $\g_n^*\to\g_{n-4}^*$
(we assume that $G_{n-4}$ is embedded into $G_n$ like
$G_{n-4}\subset G_{n-2}\subset G_n$). In particular,
$\chi_{\pi^2(f)}$ is a character of the principal series of the
group $G_{n-4}$.

On the other hand, \ref{hgh}) shows that
$\chi_{\pi^2(f)}(\pi^2(g))=\chi_{\pi^2(f)}(\pi^2(h^{-1}gh))$ for all
$h\in\w H_n$. Substituting this to (\ref{first_sreg_char_formula})
and using (\ref{second_sreg_char_formula}), we obtain
\begin{equation}
\begin{split}
&\chi(g)=q^s\cdot\chi_{\pi^2(f)}(\pi^2(g))\times\\
&\times\sum_{h\in\w H_n}\theta_{n-1}(p(h^{-1}gh))
\cdot\theta(\xi_{2,n}
\cdot(-1)^{|X_1|-1}\cdot\frac{\Delta^{X_1}(h^{-1}gh)}
{\Delta^{X_2}(h^{-1}gh)}+\xi_{n-1,n}y_{n,n-1}).
\end{split}
\label{formula_char_sreg}
\end{equation}

Furthermore,
\begin{equation}
\begin{split}&\theta_{n-1}(p(h^{-1}gh))=\theta(\xi_{1,n-1}\cdot(y_{n-1,1}
-\sum_{i=2}^{n-2}t_{n-1,i}y_{i,1})),\\
&\frac{\Delta^{X_1}(h^{-1}gh)}{\Delta^{X_2}(h^{-1}gh)}=
\frac{\Delta^{X_1}(g)+y_{n,n-1}\cdot\Delta^{X_1}(g_t)}{\Delta^{X_2}(g)},
\end{split}
\label{delta_formula}\end{equation} where $g_t$ denotes the matrix
given by replacing $y_{nj}$ by $t_{nj}$ for $2\lee j\lee n-2$ in the
matrix $g$ (the last equality can be checked directly).

Suppose now that $\chi(g)\neq0$. Thus
$\chi_{\pi^2(f)}(\pi^2(g))\neq0$ (see (\ref{formula_char_sreg})),
some conditions from the equations (\ref{sreg_uravn}), defining
$\K_D(\vfi)$, are satisfied automatically. Precisely, we have to
check only the following equalities:
\begin{equation*}
\begin{split}
&\xi_{1,n-1}y_{2,1}=\xi_{2,n}y_{n,n-1},
\quad\text{(condition (\ref{dop_usl}))}\\
&\alpha_{ij}=\beta_{ij}=0,\quad(i, j)\in D^+,\\
&\Delta_{ij}=0,\quad(i, j)\in (R(D'')\setminus
D')\cap(\w\Phi^2\cup\w\Phi_{n-1}).
\end{split}
\end{equation*}
But $(n,n-1),(2,1)\in D''$, so
$R(D'')\cap(\Phi^1\cup\Phi_n)=\varnothing$, hence, the last set of
equalities in fact has the form $\Delta_{n-1,j}=0$,
$\Delta_{i,2}=0$.

Consider the equalities $\Delta_{n-1,j}=0$ firstly. If
$j>m=\max_{(i,j)\in D}j$, then they have the form
$\Delta_{n-1,j}=y_{n-1,j}=0$, and other equalities (with $j\lee m$)
are exactly Kronecker-Capelli conditions of compatibility of the
system of linear in $t_{n-1,j}$ equations
\begin{equation}
y_{n-1,j}-\sum_{i=j+1}^{n-2}t_{n-1,i}y_{i,1}=0,\quad 2\lee j\lee m.
\label{uravn_Kron}\end{equation} In both cases equalities under
consideration are satisfied, because they express the condition
$h^{-1}gh\in P_n\rtimes\w G_{n-2}$ (see (\ref{hgh})).

On the other hand, it follows from this condition and (\ref{hgh})
that $t_{n-1,j}=-\frac{y_{n,j}}{y_{n,n-1}}$ for all $m<j\lee n-2$.
Substituting this to (\ref{uravn_Kron}), we obtain that
$\alpha_{ij}=0$ for any $(i,j)\in D^+$ (all other equalities of this
form are satisfied because of compatibility of the system
(\ref{uravn_Kron})). The equalities $\beta_{ij}=0$, $(i,j)\in D^+$,
and the other equalities $\Delta_{i,2}=0$ can be obtained by using
symmetric semi-direct decomposition of $G_n$ (see remark
\ref{nota_symm_rtimes}).

Finally, we'll prove that (\ref{dop_usl}) is satisfied. Note that
$\sum_{c\in\F_q}\theta(ct)=0$ for all $t\in\F_q^*$. Hence,  the
coefficients of $t_{n-1,j}$ in the character formula, that vary
independently, have to be zero. For example, the coefficient of
$t_{n-1,2}$ equals zero. Since $\Delta^{X_1}(g^t)=
(-1)^{|X_1|-2}\cdot t_{n-1,2}\cdot\Delta^{X_2}(g)+\ldots$ (members
without $t_{n-1,2}$), the coefficient of $t_{n-1,2}$ in the
character formula is equal to
$$-\xi_{1,n-1}y_{21}+\xi_{2,n}\cdot(-1)^{|X_1|}\cdot
y_{n,n-1}\cdot\dfrac{(-1)^{|X_1|-2}\cdot\Delta^{X_2}(g)}
{\Delta^{X_2}(g)}=-\xi_{1,n-1}y_{21}+\xi_{2,n}y_{n,n-1},$$ so,
$-\xi_{1,n-1}y_{21}+\xi_{2,n}y_{n,n-1}=0$.}\label{lemm_O_1_neq}

\lemmp{Let $f\in\mathcal O_1$ be a canonical form of a
$1$-subregular orbit, $g\in\K_D(\varphi)\subset \K_f$. Then
$\chi(g)=q^{m_D}\cdot\theta_f(e_D(\varphi))$.}{We'll compute the
value of the character under the assumption that this value isn't
equal to zero. Substituting (\ref{delta_formula}) to
(\ref{formula_char_sreg}), we have:
\begin{equation*}
\begin{split}
&\chi(g)=q^s\cdot\chi_{\pi^2(f)}(\pi^2(g))\cdot\sum_{h\in\w
H_n}\theta(\xi_{1,n-1}\cdot(y_{n-1,1}
-\sum_{i=2}^{n-2}t_{n-1,i}y_{i,1})+\\
&+\xi_{2,n} \cdot(-1)^{|X_1|}\cdot\frac{\Delta^{X_1}(g)+
y_{n,n-1}\cdot\Delta^{X_1}(g_t)}{\Delta^{X_2}(g)}+\xi_{n-1,n}y_{n,n-1}),
\end{split}
\end{equation*}
if $\chi_{\w\pi(f)}(\w\pi(g))\neq0$. If $\chi(g)\neq0$, then the
following four groups of conditions are satisfied.

1. The first condition $\chi_{\pi^2(f)}(\pi^2(g))\neq0$ is satisfied
automatically, because $g\in\K_D(\varphi)$. Indeed, it's evident,
that $\pi^2(D)\in\Phi(n-4)$ is a regular subset and
$\pi^2(g)\in\K_{\pi^2(D)}(\varphi\mid_{\pi^2(D)})$. But
$\chi_{\pi^2(f)}(\pi^2(g))$ is a character of the principal series
of $G_{n-4}$, so, Theorem \ref{char_reg} shows that
\begin{equation*}
\chi_{\pi^2(f)}(\pi^2(g))=q^{m_{\pi^2(D)}}\cdot\theta_{\pi^2(f)}(\pi^2(e_D(\phi)))
\end{equation*}
(the map $\pi^2\colon\g_n\to\g_{n-4}$ is defined by the obvious
way).

2. The second group of conditions provides the inequality
$\chi_{\w\pi(f)}(\w\pi(g))\neq0$. Since $\chi_{\pi^2(f)}(\pi^2(g))$
$\neq0$, it's enough to prove that $\Delta_{ij}(h^{-1}gh)=0$ for
$(i,j)\in (R(D'')\setminus D)\cap(\Phi^2\cup\Phi_n)$. More
precisely, these conditions have the form
\begin{equation*}
\begin{split}
&\Delta_{i,2}(g)=0,\quad (i,2)\in R(D)''\setminus D\\
&\Delta_{nj}(h^{-1}gh)=\Delta_{nj}(g)+y_{n,n-1}
\cdot\Delta_{nj}(g_t)=0,\quad (n,j)\in R(D)''\setminus D.\\
\end{split}
\end{equation*}

The first set of equalilies is contained in the Kronecker-Capelli
conditions of compatibility of the system (\ref{uravn_Kron}). Since
$t_{n-1,j}$, $3\lee j\lee m$, contain only in the second set of
equalities, it follows that $t_{n-1,j}$, $3\lee j\lee m$,
$(n-1,j)\notin R(D'')$, may be arbitrary, and other $t_{n-1,j}$ are
uniquely determined by the formulas
\begin{equation}
t_{n-1,j}=(-1)^{|X_1|-3}\cdot\frac{-\Delta_{n,j}(g)-\sum_{\substack{r>j\\(n-1,r)\notin
R(D'')}}\pm
t_{n-1,r}\cdot\Delta^{Y_r(j)}(g)}{y_{n,n-1}\cdot\Delta^{Y(j)}},
\label{small_t}\end{equation} where $Y(j)=D^+\setminus\w\Phi^j$,
$Y_r(j)=(D^+\setminus(\Phi^r\cup\w\Phi^{j-1}))\cup\{(n-r+1,j)\}$,
and signs of $t_{n-1,r}$ alternates.

3. The third condition says that we must consider $h\in\w H_n$, such
that $h^{-1}gh\in P_n\rtimes\w G_{n-2}$. This condition is always
satisfied, if $g\in\K_D(\varphi)$ (see the proof of Lemma
\ref{lemm_O_1_neq}).

4. Finally, the fourth condition says that the coefficients of
$t_{n-1,j}$ in the character formula, that vary independently, are
equal to zero. These are $t_{n-1,2}$ and $t_{n-1,j},3\lee j\lee m$,
$(n-1,j)\notin R(D'')$. The end of the proof of Lemma
\ref{lemm_O_1_neq} shows that the condition (\ref{dop_usl}) is
satisfied if and only if the coefficient of $t_{n-1,2}$ equals zero;
hence, we must prove only that the coefficients of all
$t_{n-1,j},3\lee j\lee m$, $(n-1,j)\notin R(D'')$, are equal to
zero. Simplifying (\ref{small_t}) under all our assumptions, one can
check that this is equivalent to the set of equalities
$\beta_{ij}=0$, $(i,j)\in D^+$, which are satisfied, because
$g\in\K_D(\varphi)$.

So, the character formula contains only the following expression:
\begin{equation}\begin{split}
\label{coeff_1}
&\xi_{n-1,n}y_{n,n-1}+\xi_{1,n-1}\cdot(y_{n-1,1}-\sum_{i=m+1}^{n-2}y_{i,1}\cdot(-
\dfrac{y_{n-1,i}}{y_{n,n-1}}))+\xi_{2,n}\cdot\ldots=\\
&=\xi_{n-1,n}y_{n,n-1}+\xi_{1,n-1}\cdot\dfrac{\sum_{i=n-m+1}^{n-1}y_{n,i}y_{i,1}}
{y_{n,n-1}}+\xi_{1,n-1}\cdot\dfrac{\sum_{i=m+1}^{n-m}y_{n,i}y_{i,1}}
{y_{n,n-1}}+\xi_{2,n}\cdot\ldots,\\
\end{split}
\end{equation}
where the numerator of every summand of the group marked by dots
contains the only one of the elements $y_{i,1}$, $2\lee i\lee n-m$,
and these elements aren't contained in other summands. On the other
hand, the symmetric decomposition of the group $G_n$ (see remark
\ref{nota_symm_rtimes}) give the following expression in the
character formula
\begin{equation}\begin{split}
\label{coeff_2}
&\xi_{n-1,n}y_{n,n-1}+\xi_{2,n}\cdot\dfrac{\sum_{j=2}^{n-m}y_{n,j}y_{j,1}}
{y_{21}}+\xi_{1,n-1}\cdot\ldots=\\
&=\xi_{n-1,n}y_{n,n-1}+\xi_{2,n}\cdot\dfrac{\sum_{j=2}^{m}y_{n,j}y_{j,1}}
{y_{21}}+\xi_{2,n}\cdot\dfrac{\sum_{j=m+1}^{n-m}y_{n,j}y_{j,1}}
{y_{21}}+\xi_{1,n-1}\cdot\ldots,\end{split}\end{equation} which has
to coincide with \eqref{coeff_1} (here the numerator of every
summand of the group marked by dots contains the only one of the
elements $y_{n,j}$, $n-m+1\lee j\lee n-1$, and these elements aren't
contained in other summands). Since
$\xi_{1,n-1}y_{21}=\xi_{2,n}y_{n,n-1}$, we have that, in fact,
\eqref{coeff_1} and \eqref{coeff_2} have the common form
\begin{equation*}\xi_{1,n-1}\cdot\dfrac{\sum_{i=n-m+1}^{n-1}y_{n,i}y_{i,1}}
{y_{n,n-1}}+\xi_{2,n}\cdot\dfrac{\sum_{j=2}^{n-m}y_{n,j}y_{j,1}}
{y_{21}}+\xi_{n-1,n}y_{n,n-1}
=\xi_{1,n-1}\cdot\dfrac{\gamma(g)}{y_{n,n-1}}+\xi_{n-1,n}y_{n,n-1}.
\end{equation*}

We also have
\begin{equation}
\label{q_deg}
\sum_{h\in\w
H_n}\theta(\xi_{1,n-1}\cdot\dfrac{\gamma(g)}{y_{n,n-1}}+\xi_{n-1,n}y_{n,n-1})=
q^{s_1}\cdot\theta(\xi_{1,n-1}\cdot\dfrac{\gamma(g)}{y_{n,n-1}}+\xi_{n-1,n}y_{n,n-1}),
\end{equation}
where $s_1=|R(D)\cap\Phi_{n-1}\cap\w\Phi^m|$. Indeed, $t_{n-1,j}$,
$2\lee j\lee m$, $(n-1, j)\in R(D)$, vary independently, and if we
fix them, then all other $t_{n-1,j}$ are uniquely determined.

Note that
\begin{equation}
\label{coeff_raven_eDphi}
\begin{split}
&\xi_{1,n-1}\cdot\dfrac{\gamma(g)}{y_{n,n-1}}=
\xi_{1,n-1}\cdot\dfrac{\gamma(x_D(\varphi))}{y_{n,n-1}}=
\xi_{1,n-1}\cdot\dfrac{a_{n-d+1,n-d}a_{n-d,1}+
a_{n-d+1,2}a_{2,1}}{a_{n,n-1}}=\\
&=\xi_{1,n-1}a_{n-1,1}+ \dfrac{\xi_{1,n-1}}{a_{n,n-1}}\cdot
a_{n-d+1,2}a_{2,1}= \xi_{1,n-1}a_{n-1,1}+
\dfrac{\xi_{2,n-d+1}}{a_{2,1}}\cdot a_{n-d+1,2}a_{2,1}=\\
&=\xi_{1,n-1}a_{n-1,1}+\xi_{2,n-d+1}a_{n-d+1,2}
\end{split}
\end{equation}
(here $a_{ij}=\varphi(i, j)$, as above).

Finally, substitute (\ref{q_deg}) and (\ref{coeff_raven_eDphi}) to
(\ref{formula_char_sreg}):
\begin{equation*}
\begin{split}
&\chi_{n,f}(g)=q^s\cdot\chi_{\pi^2(f)}(\pi^2(g))\cdot\sum_{h\in\w
H_n}\theta(\xi_{1,n-1}\cdot\dfrac{\gamma(g)}{y_{n,n-1}}+\xi_{n-1,n}y_{n,n-1})=\\
&=q^{s+s_1}\cdot\chi_{\pi^2(f)}(\pi^2(g))\cdot\theta(\xi_{n-1,1}a_{n,n-1}+
\xi_{1,n-1}a_{n-1,1}+\xi_{2,n-d+1}a_{n-d+1,2})=q^{m_D}\cdot\theta_f(e_D(\varphi)),
\end{split}
\end{equation*}
because $m_{\pi^2(D)}+s+s_1=m_D$. The proof is complete.}

\lemmp{Let $1<d<n_1$, $f\in\mathcal O_d$ be a canonical form of a
$d$-subregular, $\chi=\chi_f$ be the corresponding character, and
$g\in G_n$. Then $\chi(g)\neq0$ if and only if
$g\in\K_D(\varphi)\subset\K_f$; in this case
$\chi(g)=q^{m_D}\cdot\theta_f(e_D(\varphi))$.}{Suppose that the
conditions of the Lemma are satisfied. According to
(\ref{ind_reps}),
\begin{equation*}
T_f=\Ind{P_n\rtimes G_{n-2}}{G_n}{\theta_n\otimes T_{\pi(f)}},
\end{equation*}
hence,
\begin{equation*}
\chi(g)=\sum_{h\in
H_n}\theta_n(p(h^{-1}gh))\cdot\chi_{\pi(f)}(\pi(h^{-1}gh))
\end{equation*}
(here the summation is over $h\in H_n$ such that $h^{-1}gh\in
P_n\rtimes G_{n-2}$). Since $\pi(f)$ is a canonical form of a
subregular orbit of $G_{n-2}$, we may assume that Theorem
\ref{theo_sreg} is valid for $\chi_{\pi(f)}$ (inductive assumption).

Using the same arguments, as in two previous Lemmas
(Kronecker-Capelli criterion of compatibility of systems of linear
equations, conditions of the form $h^{-1}gh\in P_n\rtimes G_{n-2}$,
and vanishing the coefficients of independent variables in the
character formula), one can obtain the required
result.}\label{lemm_for_all_d}

\sect{Remarks and generalizations}\label{utochn_obob} Above we
considere $d$-subregular subsets $D$ with $D(d)=D_1(d)$ and assume
that $d<n_1$. Here we'll consider all other cases.

First, if $D\supset D_0(d)$, then the subset $D$ is basic and the
conjugacy class $\K_D(\varphi)$ is given by the equations
(\ref{reg_uravn}) (see \cite{Andre2}).

Second, we should consider subregular orbits of the second and the
third types (these are $n_1$-subregular orbits for even and odd $n$
respectively). Note that if $f=(\xi_{ij})\in\mathcal O_{n_1}$ is a
canonical form of a subregular orbit of the second type, then the
set $D=\{(i,j)\in\Phi(n)\mid\xi_{ji}\neq0)\}$ is basic itself, so in
this case the subregular character is a basic character (in the
sense of Andre); a complete description of such characters is given
in \cite{Andre2}. This is also true for subregular orbits of the
third type, if $\xi_{n_1,n_0+1}=0$ (see definition
\ref{defi_sreg_orbits}).

If $\Ou_f$ is a subregular orbit of the third type and
$\xi_{n_1,n_0+1}\neq0$, then all statements are the same as in the
case of orbits of the first type (except that the ideal $J$ doesn't
contain the polynomial $\gamma-c_0$). All proofs are quite similar
to proofs in case of orbits of the first type, so we don't reply
them.

Finally, we should explain the genesis of polynomials $\alpha_{ij}$,
$\beta_{ij}$ and $\gamma$ from our description of subregular
characters of the first type with $D\supset D_1(d)$: (these
polynomials are \emph{not} minors!). Consider the following
\emph{characteristic matrix}
\begin{equation*}
M(t)=\begin{pmatrix}
1&0&\ldots&0&0\\
t\cdot y_{21}&1&\ldots&0&0\\
\vdots&\vdots&\ddots&\vdots&\vdots\\
t\cdot y_{n-1, 1}&t\cdot y_{n-1, 2}&\ldots&1&0\\
t\cdot y_{n, 1}&t\cdot y_{n, 2}&\ldots&t\cdot y_{n, n-1}&1
\end{pmatrix}.
\end{equation*}

A minor of the matrix $M(t)$ is polynomial on $t$, and its
coefficients are polynomial on $y_{ij}$ (hence, we can consider them
as function on $G_n$ or as elements of $\F_q[\g_n]$). \deno{By
$M^{i_1,\ldots,i_k}_{j_1,\ldots,j_k}(t^d, g)$ we denote the value of
the coefficient of $t^d$ in the minor of the matrix $M(t)$ with rows
$1\lee i_1,\ldots,i_k\lee n$ and columns $1\lee j_1,\ldots,j_k\lee
n$ on an element $g=(y_{ij})\in G_n$. We also denote
$M^X(d,g)=M^{\sigma(i_1,\ldots,i_k)}_{\tau(j_1,\ldots,j_k)}(t^d,
g)$, where $X$ is the set of pairs of the form
$X=\{(i_1,j_1),\ldots,(i_k,j_k)\}$, and $\sigma, \tau$ are
permutations such that $\sigma(i_1)<\ldots<\sigma(i_k)$ and
$\tau(j_1)<\ldots<\tau(j_k)$.}

For example, for any $1\lee d\lee n_0$ we have $M^X(d,
g)=\Delta_d(g)$, where $X=\{(n, 1), (n-1, 2),\ldots,(n_1+2,n_0)\}$;
hence, all minors considered above are also coefficients of minors
of characteristic matrix. On the other hand, if $(i,j)\in D^+$
(where $D\supset D_1(d)$ is a $d$-subregular subset and $1\lee
d<n_1$), $m=\max_{(i,j)\in D}j$, then
\begin{equation*}
\begin{split}
&\alpha_{ij}(g)=\pm M^{X^{\alpha}_{ij}}(2,g), \quad
X^{\alpha}_{ij}=\{n-d+1,j\}\cup\bigcup_{d<i\lee m}\{(i,i)\},\\
&\beta_{ij}(g)=\pm M^{X^{\beta}_{ij}}(2,g), \quad
X^{\beta}_{ij}=\{i,d\}\cup\bigcup_{d<j\lee m}\{(j,j)\},\\
&\gamma(g)=\pm M^{X^{\gamma}}(2,g), \quad
X^{\gamma}=\{n-d+1,d\}\cup\bigcup_{d<i<n-d+1}\{(i,i)\}\\
\end{split}
\end{equation*}
(the choice of signs depends on the cardinality of $D^+$).

We see that subregular characters can be described in terms of
coefficients of minors of the characteristic matrix. In fact,
subregular orbits, all orbits for $n\lee7$ \cite{IgnatevPanov} and
all irreducible characters for $n\lee5$ \cite{Ignatev} can be
described in these terms too. So, we may conjecture that all orbits
and characters of the unitriangular group of an arbitrary dimension
can be describes in terms of coefficients of minors of the
characteristic matrix.

The author expresses his sincere gratitude to professor A.N.~Panov
for suggesting the ideas that gave rise to this work.

Department of Algebra and Geometry, Samara State University, ul.
Akad. Pavlova, 1, Samara 443011, Russia

E-mail address: \texttt{mihail.ignatev@gmail.com}

\end{document}